\documentclass[12pt,reqno]{amsart}
\usepackage{amsfonts,amssymb,amsthm,enumerate,color,comment, pgf,tikz}
\usepackage[left=2.6cm,right=2.6cm,top=3cm,bottom=3cm,bindingoffset=0cm]{geometry}
\newtheorem{thm}{Theorem}[section]
\newtheorem{lem}[thm]{Lemma}

\newtheorem{prop}[thm]{Proposition}
\newtheorem*{claim}{Claim}
\newtheorem*{claim1}{Claim 1}
\newtheorem*{claim2}{Claim 2}
\newtheorem*{claim3}{Claim 3}
\newtheorem*{claim4}{Claim 4}

\newtheorem{hypo}[thm]{Hypothesis}
\theoremstyle{remark}
\newtheorem{rem}[thm]{Remark}
\theoremstyle{definition}
\newtheorem{defi}[thm]{Definition}

\def\A{\mathcal{A}}
\def\B{\mathcal{B}}

\def\G{\Gamma}
\def\k{\mathbf{k}}
\def\K{\mathbf{K}}
\def\S{\mathcal{S}}
\def\Z{\mathbb{Z}}
\DeclareMathOperator{\aut}{Aut}
\DeclareMathOperator{\cay}{Cay}
\DeclareMathOperator{\orb}{Orb}
\DeclareMathOperator{\Span}{Span}
\DeclareMathOperator{\rk}{rank}
\DeclareMathOperator{\rad}{rad}
\DeclareMathOperator{\Sup}{Sup}
\DeclareMathOperator{\sym}{Sym}
\DeclareMathOperator{\vertex}{V}
\DeclareMathOperator{\edge}{E}
\newcommand{\sg}[1]{\langle {#1}\rangle}
\begin{document}
\title[Stability of Cayley graphs and Schur rings]
{Stability of Cayley graphs and Schur rings}
\author[A. Hujdurovi\'c \and I. Kov\'acs]
{Ademir Hujdurovi\'c$^{\, 1}$ \and Istv\'an Kov\'acs$^{\, 2}$}
\address{A. Hujdurovi\'c \and I.~Kov\'acs 
\newline\indent
UP IAM, University of Primorska, Muzejski trg 2, SI-6000 Koper, Slovenia 
\newline\indent
UP FAMNIT, University of Primorska, Glagol\v jaska ulica 8, SI-6000 Koper, Slovenia}
\email{ademir.hujdurovic@upr.si}
\email{istvan.kovacs@upr.si} 
% support
\thanks{$^1$~Supported by the Slovenian Research Agency (research program P1-0404, research projects N1-0140, N1-0102, J1-1691, N1-0159, J1-2451, N1-0208, J1-4084). 
\newline\indent 
$^2$~Supported by the Slovenian Research Agency 
(research program P1-0285 and research projects J1-1695, N1-0140, J1-2451, N1-0208 and J1-3001).}
%--------------------------------------------------------------------------------------------------------------%
\maketitle
\begin{abstract}
A graph $\G$ is said to be unstable if for the direct product $\G \times K_2$, 
$\aut(\G \times K_2)$ is not isomorphic to $\aut(\G) \times  \Z_2$. 
In this paper we show that a connected and non-bipartite Cayley graph $\cay(H,S)$ is unstable if and only if the set $S \times \{1\}$ belongs to a Schur ring over the group $H \times \Z_2$ having certain properties. The Schur rings with these properties are characterized if $H$ is an abelian group of odd order or a cyclic group of twice odd order. As an application, a short proof is given 
for the result of Witte Morris stating that every connected unstable Cayley graph on an abelian group of odd order has twins (Electron.~J.~Combin, 2021). 
As another application, sufficient and necessary conditions are given for a 
connected and non-bipartite circulant graph of order $2p^e$ to be unstable, where 
$p$ is an odd prime and $e \ge 1$. 
\end{abstract}
%--------------------------------------------------------------------------------------------------------------%
\section{Introduction}\label{sec:intro}

All groups in this paper will be finite and all graphs will be finite and simple.  
If $\G$ is a graph, then $\vertex(\G)$, $\edge(\G)$ and $\aut(\G)$ denote its vertex set, edge set and automorphism group, 
respectively. The direct product $\G \times K_2$ of a graph $\G$ and the complete 
graph $K_2$ on two vertices, also known as the \emph{canonical double cover} of 
$\G$, is defined to have vertex set $\vertex(\G) \times \{0,1\}$ and edges 
$\{(u,0),(v,1)\}$, where $\{u,v\} \in \edge(\G)$. The graph $\G \times K_2$ admits natural automorphisms, namely the permutation 
$$
(v,i) \mapsto (v,1-i),~\text{where}~v \in \vertex(\G),~i=0,1;
$$ 
and for every $\alpha \in \aut(\G)$, the permutation 
$$
(v,i) \mapsto (v^\alpha,i),~\text{where}~v \in \vertex(\G),~i=0,1.
$$
These permutations can be easily seen to form a group, which is isomorphic to 
$\aut(\G) \times \Z_2$. Now, we say that $\G$ is \emph{stable} if 
$\aut(\G \times K_2) \cong \aut(\G) \times \Z_2$, and \emph{unstable} otherwise. 
This concept of stability was defined in \cite{MSZ}. 
Recently, several papers were devoted to the stability of graphs~\cite{HuMi,QXZ21,
W-M23}, especially to circulant graphs~\cite{FH,HMW21a,HMW21b,QXZ19,W}.  

Let $H$ be a group with identity element $1_H$, and let $S \subseteq H$ be a subset 
such that $1_H \notin S$ and $x^{-1} \in S$ whenever $x \in S$. The \emph{Cayley graph} $\cay(H,S)$ is defined to have vertex set $H$ and edges 
$\{x,sx\}$, where $x \in H$ and $s \in S$. In the case when $H$ is a cyclic group, the term \emph{circulant graph} is commonly used.  

The goal of this paper is to propose an approach to the stability of Cayley 
graphs using Schur rings or S-rings for short. S-rings over a group $G$ are 
certain subrings of the integer group ring $\Z G$, which were defined by 
Wielandt~\cite{Wbook}, and first studied by Schur in his investigation of permutation 
groups. S-rings also serve as an effective tool in algebraic combinatorics~\cite{MP}. For the exact definition of an S-ring and the definitions of all the S-ring theoretical concepts, which appear in our main results below, we refer to Section~\ref{sec:S-rings}. 
\medskip

In Section~\ref{sec:main1}, we establish a sufficient and 
necessary condition for a connected and non-bipartite Cayley graph to be unstable in 
terms of S-rings. We remark that none of the latter constrains on the given Cayley 
graph is essential. It is easy to show that any disconnected graph as well as any bipartite graph with a non-trivial automorphism group is unstable (see, e.g.,~\cite{W}). 

\begin{thm}\label{main1}
Let $G=H \times \sg{a}$, where $H$ is any finite group and $\sg{a} \cong \Z_2$. 
The following conditions are equivalent for every connected non-bipartite graph 
$\cay(H,S)$, 
\begin{enumerate}[(1)]
\item $\cay(H,S)$ is unstable. 
\item There exists a Schurian S-ring $\A$ over $G$ such that 
$\underline{H}, \underline{Sa} \in \A$ and $\underline{\{a\}} \notin \A$.
\end{enumerate}
\end{thm}

Theorem~\ref{main1} suggests the following recipe for finding all connected and non-bipartite unstable graphs $\cay(H,S)$ on a given group $H$.
\medskip

\begin{tabular}{lp{12cm}}
{\bf Step~1.} & Let $G=H \times \sg{a}$, where $\sg{a} \cong \Z_2$. Describe   
all S-rings $\A$ over $G$ with $\underline{H} \in \A$ and $\underline{a} \notin \A$. 
\\ [+1ex]
{\bf Step~2.} & Describe the connected non-bipartite graphs $\cay(H,S)$ using the 
fact that  $\underline{Sa} \in \A$ for some S-ring $\A$ described in Step~1. 
\end{tabular}
\medskip

%By Step~1 one can take advantage of the theory of S-rings, which has been especially %developed for abelian groups. 
In this paper we explore this idea in two cases: 
$H$ is an abelian group of odd order and $H$ is a cyclic group of order $2p^e$, where 
$p$ is an odd prime and $e \ge 1$.   
\medskip

In Section~\ref{sec:main2+W-M}, we focus on the case when $H$ is an abelian group of odd order. The following theorem will be derived about the S-rings described in Step~1. 

\begin{thm}\label{main2}
Let $G=H \times \sg{a}$, where $H$ is an abelian group of odd order and 
$\sg{a} \cong \Z_2$. If $\A$ is an S-ring over $G$ with $\underline{H} \in \A$ and $\underline{a} \notin \A$, then $\A$ is the $H/L$-wreath product for some subgroup 
$L \le H$, $L \ne 1$.
\end{thm}

Theorem~\ref{main2} will follow directly from a result of Somlai and 
Muzychuk~\cite{SM} on S-rings over abelian groups of twice odd order. 
Then in accordance with Step~2, we obtain a short 
proof for the following recent result of Witte Morris~\cite{W-M21}. 

\begin{thm}[\cite{W-M21}]\label{W-M}
If $H$ is an abelian group of odd order, then every unstable connected graph 
$\cay(H,S)$ has twins. 
 \end{thm}

Two vertices of a graph form \emph{twins} if they are adjacent with the same vertices. 
We would like to emphasise that the proof of Theorem~\ref{W-M} given by Witte Morris is entirely different from our proof.
\medskip

Section~\ref{sec:main3} is devoted to the S-rings mentioned in Step~1 in the case when 
$H$ is a cyclic group of twice odd order. The main result of the section is 
the following theorem.

\begin{thm}\label{main3}
Let $G=H \times \sg{a}$, where $H \cong \Z_{2n}$, $n>1$, $n$ is odd, and 
$\sg{a} \cong \Z_2$. If $\A$ is an S-ring over $G$ with 
$\underline{\sg{H}} \in \A$ and $\underline{\{a\}} \notin \A$, then 
$\{a,ab\}$ is a basic set of $\A$, or 
$$
\bigcap_{X \in \S(\A) \atop X \cap H_0a \ne \emptyset}\rad(X \cap H_0a) \ne 1,
$$
where $b$ is the unique involution of $H$ and $H_0$ is the subgroup of $H$ of order 
$n$. 
\end{thm}

Finally, we are going to use Theorem~\ref{main3} to study the stablility of circulant graphs. 
To the best of our knowledge, a characterisation of the unstable circulant 
graphs of order $n$ and valency $k$ is known only in special cases:  
$n$ is odd~\cite{FH,QXZ19}, $n=2p$ for a prime $p$~\cite{HMW21a}, or 
$k \le 7$~\cite{HMW21b}. Here we handle the case when the order 
$n=2p^e$ for an odd prime $p$ and $e \ge 1$.

\begin{thm}\label{main4}
Let $H \cong \Z_{2p^e}$, where $p$ is an odd prime and $e \ge 1$. 
A connected and non-bipartite 
graph $\cay(H,S)$ is unstable if and only if 
\begin{enumerate}[(1)]
\item $e > 1$ and $(S \cap H_0)h=S \cap H_0$, where 
$H_0$ is the unique subgroup of $H$ of order $p^e$ and $h \in H_0, h \ne 1_G$; or  
\item $\cay(H,S) \cong \cay(H,Sb)$, where $b$ is the unique involution of $H$.
\end{enumerate}
\end{thm}

\begin{rem}\label{to:main4}
In fact, the sufficiency part of the theorem follows from known constructions of 
unstable circulant graphs. The graphs satisfying the condition in (1) are of 
Wilson type (C.1) (see~\cite{W}),
and those satisfying the condition in (2)  
are unstable by \cite[Proposition~3.7]{HMW21a}. 
\end{rem}

Besides S-ring theory, another ingredient in our proof of Theorem~\ref{main4} is Muzychuk's remarkable solution for the isomorphism problem of circulant graphs~\cite{M}.  
A short account on this result can be found in Section~\ref{sec:iso}. The proof of 
Theorem~\ref{main4} will be given in Section~\ref{sec:main4}.
%--------------------------------------------------------------------------------------------------------------%
\section{Schur rings}\label{sec:S-rings}

Let $G$ be a group with identity element $1_G$. 
The set of the non-identity elements will be denoted by $G^\#$ and 
for $g \in G$, let $o(g)$ denote the order of $g$. 
For a permutation group $A \le \sym(G)$ and $g \in G$, denote by $A_g$ the \emph{stabiliser} of $g$ in $A$ and by $\orb_A(g)$ the \emph{orbit} of $g$ under $A$. 
We let $\orb(A,G)=\{\orb_A(x): x \in G\}$.
For $g \in G$, the \emph{right multiplication} $g_r$ is 
the permutation of $G$ acting as $x \mapsto xg$ ($x \in G$). We let 
$G_r=\{ x_r : x \in G\}$ and set 
$$
\Sup(G)=\{ A \le \sym(G) : G_r \le A\}.
$$

For a subset $X \subseteq G$, the element $\sum_{x \in X}x$ in the group 
ring $\Z G$ is denoted by $\underline{X}$. 

\begin{defi}
(Wielandt~\cite[Chapter~IV]{Wbook}) 
A subring $\A$ of the integer group ring $\Z G$ is called a 
\emph{Schur ring} (\emph{S-ring} for short) if there exists a partition
$\S(\A)$ of $G$ such that 

\begin{enumerate}[(1)]
\setlength{\itemsep}{0.4\baselineskip}
\item $\{1_G\} \in \S(\A)$.
\item  If $X \in \S(\A)$ then $X^{-1} \in \S(\A)$.
\item $\A=\Span_{\Z}\{\underline{X} :\ X \in \S(\A)\}$.
\end{enumerate}
\end{defi}

The subsets in $\S(\A)$ are called the \emph{basic sets} of
$\A$ and the number $\rk(\A):=|\S(\A)|$ is called the \emph{rank} of $\A$. 
The motivation of the above definition can be explained by the following result of 
Schur. 

\begin{thm}[\cite{Sch}]
If $A \in \Sup(G)$, then the free $\Z$-module
$\Span_{\Z}\{\underline{X} : X \in \orb(A_{1_G},G)\}$ is a subring of $\Z G$.
\end{thm}

The ring in the theorem is an example of an S-ring, which is also 
called the \emph{transitivity module} over $G$ induced by $A$ and denoted by
$V(G,A_{1_G})$.  An S-ring $\A$ is called \emph{Schurian} if $\A=V(G,B_{1_G})$ for some permutation group $B \in \Sup(G)$.   
We remark that not all S-rings are Schurian (see \cite{Wbook}).

If $\A$ and $\B$ are two S-rings over $G$, then their usual 
intersection $\A \cap \B$ is also 
an S-ring over $G$ (see, e.g.,~\cite[the paragraph following Theorem~4.2]{M}).  Moreover, if both $\A$ are $\B$ are Schurian, 
then $\A \cap \B$ is also Schurian. 
\medskip

Let $\A$ be an S-ring over a group $G$. 
A subset $X \subseteq G$ is called an
\emph{$\A$-set} if $\underline{X} \in \A$, and a subgroup
$H \le G$ is called an \emph{$\A$-subgroup} if
$\underline{H} \in \A$. We say that $\A$ is \emph{primitive} if $1$ and $H$ are the only $\A$-subgroups of $H$. 

There are two natural $\A$-subgroups associated with an $\A$-set $X$, namely, 
$\sg{X}$ and the \emph{radical} of $X$ defined as 
$$
\rad(X)=\{g \in G : Xg=X~\text{and}~gX=X\}
$$
(see \cite[Propositions~23.6~and~23.5]{Wbook}). 
If $H$ and $K$ are two $\A$-subgroups, then it can be easily checked that so are 
$H \cap K$ and $\sg{H \cup K}$. 

Let $H \le G$ be an $\A$-subgroup. 
Then the free $\Z$-module 
$$
\A_H:=\Span_{\Z}\{ \underline{X} : X \in \S(\A), X \subseteq H \}
$$ 
is an S-ring over $H$, which is called an \emph{induced S-subring} of $\A$. 
Furthermore, if $X \in \S(\A)$, then there is a positive constant $\ell$ such that 
\begin{equation}\label{eq:cap}
\forall x \in G :~|Hx \cap X|=0~\text{or}~\ell. 
\end{equation}
%(for a proof, see \cite[pp.~167--168]{FK}).

Assume, in addition, that $H \trianglelefteq G$. 
For an arbitrary subset $X \subseteq G$, we let $X/H$ denote the subset of the 
quotient group $G/H$ defined by 
$$
X/H=\{Hx : x \in X\}.
$$  
It follows that the sets $X/H$ form the basic sets of an S-ring over $G/H$ while $X$ runs 
over $\S(\A)$ (see~\cite{T}). The latter S-ring is called \emph{a quotient S-ring} and 
denoted by $\A_{G/H}$. 
In what follows, if $K, L$ are two $\A$-subgroups such that $L \trianglelefteq K$, 
then the more simple notation $\A_{K/L}$ will be used instead of $(\A_K)_{K/L}$. 
Note that, if $\A$ is Schurian, then so is $\A_{K/L}$.  
\medskip 

For a subset $X \subseteq G$ and integer $m$, define
$X^{(m)}=\{ x^m : x \in X\}$, and for a group ring
element $\eta=\sum_{x \in G}c_x x$, define
$\eta^{(m)}=\sum_{x \in G}c_x x^m$. 
If $G$ is also abelian and $d$ is a divisor of $|G|$, then we use 
the notation $G[d]=\{x \in G : o(x)=d\}$, and if $d=p$ is a prime we define  
$$
X^{[p]}:=\{ x^p : x \in X~\text{and}~|X \cap xG[p]| \not\equiv 0\!\!\!\!\pmod p\}.
$$
Note that $X^{[p]}$ is possibly the empty set. 
The next properties, also referred to as Schur's first and second theorem 
on multipliers, respectively (see~\cite{MP}), will be used often in the sequel. 

\begin{thm}[{\cite[Theorem~23.9]{Wbook}}]\label{S-multi}
Let $\A$ be an S-ring over an abelian group $G$.
\begin{enumerate}[(1)]
\item If $m$ is an integer coprime to $|G|$ and $\eta \in \A$, then
$\eta^{(m)} \in \A$. In particular, $X^{(m)} \in \S(\A)$ whenever $X \in \S(\A)$.
\item If $p$ is a prime divisor of $|G|$ and $X$ is an $\A$-set, 
then $X^{[p]}$ is an $\A$-set. 
\end{enumerate}
\end{thm}

We proceed with the star and the generalised wreath product of S-rings. 
The former was introduced by Hirasaka and Muzychuk~\cite{HiMu} and the latter by
Evdokimov and Ponomarenko~\cite{EP02} and independently by
Leung and Man~\cite{LM} under the name  wedge product.

\begin{defi}[\cite{HiMu}]\label{star}
Let $\A$ be an S-ring over a group $G$ and
$V, W \leq G$ be two $\A$-subgroups. The S-ring
$\A$ is the \emph{star product} of $\A_V$ with $\A_W$,
written as $\A=\A_V \star \A_W$, if
\begin{enumerate}[(1)]
\setlength{\itemsep}{0.4\baselineskip}
\item $V \cap W \trianglelefteq W$.
\item Every $X \in \S(\A)$, $X \subseteq W \setminus V$ is a union of some
$(V \cap W)$-cosets.
\item For every $X \in \S(\A)$ with $X \subseteq G\setminus
(V \cup W)$, there exist basic sets $Y, Z \in \S(\A)$
such that $X=Y Z$, $Y \subseteq V$ and $Z \subseteq W$.
\end{enumerate}
\end{defi}

The star product is \emph{non-trivial} if
$1 < V <  G$.  In the special case when
$V \cap W=1$ it is also called the \emph{tensor product} and written
as $\A_V \otimes \A_W$.
\medskip

If $\A$ and $\B$ are two S-rings over $G$ such that $\A \subseteq \B$, 
then $\A$ is also called an \emph{S-subring} of $\B$. In this case every basic set of $\A$ 
can be written as a union of basic sets of $\B$.   

\begin{defi}[\cite{EP02}]\label{gwp}
Let $\A$ be an S-ring over a group $G$ and let $L, U$ be 
$\A$-subgroups of $G$ such that $L \le U$.  
The S-ring $\A$ is the \emph{$U/L$-wreath product} (also called the \emph{generalised wreath product} of
$\A_U$ with $\A_{G/L}$) %, written as $\A=\A_U \wr_{U/L} \A_{G/L}$,  
if
\begin{enumerate}[(1)]
\setlength{\itemsep}{0.4\baselineskip}
\item $L \trianglelefteq G$.
\item For every $X \in \S(\A)$, $X \subseteq G \setminus U$ is a union of some $L$-cosets.
\end{enumerate}
\end{defi}

The $U/L$-wreath product is \emph{non-trivial} if $L \neq 1$ and $U \neq G$. 
%In the special case when $U=L$, the $U/L$-wreath product is also called \emph{wreath product} and written as $\A_U \wr \A_{G/U}$.
The following simple relation with the star product will be used later, hence 
we record it here.  

\begin{lem}\label{star-is-gwp}
Let $\A$ be an S-ring over a group $G$ such that 
$\A=\A_V \star \A_W$ and $V \cap W \trianglelefteq G$. 
Then $\A$ is the $V/(V \cap W)$-wreath product.
\end{lem}
\begin{proof}
Let $X \in \S(\A)$ be an arbitrary basic set outside $V$. 
We have to show that $V \cap W \le \rad(X)$. 
This follows from Definition~\ref{star}(2) if $X \subseteq W$. 
Let $X$ be outside $W$. By Definition~\ref{star}(3), 
there exist basic sets $Y, Z \in \S(\A)$
such that $X=Y Z$, $Y \subseteq V$ and $Z \subseteq W$.
Then $V \cap W \le \rad(Z)$, implying that $V \cap W \le \rad(YZ)=\rad(X)$.  
\end{proof}
 
In the rest of the section we review the structure theorems of S-rings over 
cyclic $p$-groups obtained by P\"oschel~\cite{P}, where 
$p$ is an odd prime. 
For an integer $e \ge 1$, let $[e]=\{1,\ldots,e\}$. 
We say that a partition $\Pi$ of $[e]$ is an \emph{interval partition} if every class of 
$\Pi$ is in the form $\{i,i+1\ldots,j\}$, where $1 \le i \le j \le e$.
Given a positive integer $n$, denote by $n_p$ and $n_{p'}$ the 
\emph{$p$-part} and the \emph{$p'$-part} of $n$, respectively, i.e., 
$n_p$ is the largest power of $p$ that divides $n$, and $n_{p'}=n/n_p$. 

\begin{defi}\label{S-system}
Let $p$ be an odd prime and $e \ge 1$ be an integer.  
By an \emph{S-system} with respect to $p^e$ we mean an ordered 
$(e+1)$-tuple $(d_1,\ldots,d_e;\Pi)$, where for every $i$, $1 \le i \le e$, 
$d_i$ is a divisor of $p^{i-1}(p-1)$, $\Pi$ is an interval partition of $[e]$, and the 
following conditions hold:
\begin{enumerate}[(1)]
\setlength{\itemsep}{0.25\baselineskip}
\item If $I \in \pi$, $|I| > 1$, and $i \in I$, then $d_i=p^{i-1}(p-1)$.
\item If $2 \le i \le e$, then $(d_{i-1})_p \le (d_i)_p$.  
\item If $2 \le i \le e$ and $(d_i)_p < p^{i-1}$, then 
$(d_{i-1})_{p'}=(d_{i})_{p'}$.
\end{enumerate}
\end{defi}

Suppose that $G \cong \Z_{p^i}$ for an odd prime $p$ and integer $i \ge 1$. 
It is well-known that $\aut(G) \cong \Z_{p^{i-1}(p-1)}$, and therefore, for every divisor $d$ of $p^{i-1}(p-1)$, there is a unique subgroup of $\aut(G)$ of order $d$. 
Recall that
$
G[p^j]=\{x \in G : o(x)=p^j\},
$ for any $j\leq i$.

The next theorem contains a construction of S-rings based on S-systems. 

\begin{thm}[{\cite[(4.13)]{P}}]\label{P1}
Let $G=\sg{a}  \cong \Z_{p^e}$ for an odd prime $p$ and $e \ge 1$, and let  
$(d_1,\ldots,d_e;\Pi)$ be an S-system with respect to $p^e$. 
For the integer $i$, $1 \le i \le e$, let $K_i \le \aut(\sg{a^{p-i}})$ be the unique subgroup of order $d_i$. Then the free $\Z$-module 
$$
\Span_{\Z}\big\{ \underline{X} : X \in \Delta\big\} 
$$
is a Schurian S-ring over $G$, where $\Delta$ is the partition of $G$ defined as 
\begin{equation}\label{eq:Delta}
\Delta=\big\{\, \{1_G\}\, \big\}~\cup~
\Big\{\, \bigcup_{i \in I}G[p^i] : I \in \Pi,  |I| > 1\, \Big\}~\cup~ 
\Big\{\, \orb_{K_i}(x) : \{i\} \in \Pi, x \in G[p^i]\, \Big\}.
\end{equation}
\end{thm}

The S-ring described in Theorem~\ref{P1} will be denoted by $\S_G(d_1,\ldots,d_e;\Pi)$. 
P\"oschel~\cite{P} showed that the S-rings over cyclic $p$-groups are exactly those 
constructed from the S-systems. 

\begin{thm}[{\cite[(4.12)]{P}}]\label{P2}
Let $G \cong \Z_{p^e}$ for an odd prime $p$ and $e \ge 1$
If $\A$ is an S-ring over $G$, then $\A=\S_G(d_1,\ldots,d_e;\Pi)$ for some  
S-system $(d_1,\ldots,d_e;\Pi)$ with respect to $p^e$. 
\end{thm}
%--------------------------------------------------------------------------------------------------------------%
\section{Proof of Theorem~\ref{main1}}\label{sec:main1}

Let $G=H \times \sg{a}$ where $H$ is any group and $\sg{a} \cong \Z_2$, and 
let $\cay(H,S)$ be a connected non-bipartite graph. 
It can be easily seen that $\cay(H,S) \times K_2 \cong \cay(G,Sa)$. Moreover, 
$\cay(H,S)$ is stable if and only if 
\begin{equation}\label{eq:stab}
\aut(\cay(G,Sa))=\aut(\cay(H,S)) \times \sg{a_r},
\end{equation}
where by the latter group we mean direct product of two permutation groups acting 
on $G=H \times \sg{a}$. 
For the sake of simplicity we set $A=\aut(\cay(G,Sa))$ and write $1$ for $1_G$. 

\begin{claim} 
$\cay(H,S)$ is stable if and only if $a_r \alpha = \alpha a_r$ for every $\alpha \in A_1$.  
\end{claim}
\begin{proof}[Proof of the claim]
The implication ``$\Rightarrow$'' is clear by \eqref{eq:stab}. 

For the implication ``$\Leftarrow$'' assume that  $a_r \alpha = \alpha a_r$ for every $\alpha \in A_1$. 
The graph $\cay(G,Sa)$ is bipartite with colour classes $H$ and $Ha$. 
Since $\cay(H,S)$ is connected and non-bipartite, it follows that $\cay(G,Sa)$ is also connected. Therefore, the partition of $G$ into $H$ and $Ha$ is $A$-invariant.
Let $\alpha \in A_1$. Then $H^\alpha=H$. Let $\beta$ be the permutation of 
$H$ induced by $\alpha$. Then for every $x \in H$, $(xa)^\alpha=x^{a_r\alpha}=
x^{\alpha a_r}=x^{\beta}a$. This means that 
$\alpha \in \sym(H) \times \sg{a_r}$. 
We show below that $\beta \in \aut(\cay(H,S))$.  

Pick an arbitrary edge $\{x,sx\}$ of $\cay(H,S)$. 
Then $\{x,sax\} \in \edge(\cay(G,S))$, and since $\alpha \in A$, it follows 
that  
$$
(sax)^\alpha=s'ax^\alpha~\text{for some}~s' \in S.
$$ 

On the other hand, $(sax)^\alpha=(sx)^\beta a$ and $s'ax^\alpha=s' x^\beta a$. 
We obtain that $\beta$ maps the edge $\{x,sx\}$ to the edge 
$\{x^\beta,s' x^\beta\}$, so $\beta \in \aut(\cay(H,S))$. We showed that 
$A_1 \le \aut(\cay(H,S)) \times \sg{a_r}$. Using this, together with the fact that $A=A_1 G_r$ and \eqref{eq:stab}, we deduce that $\cay(H,S)$ is stable. 
\end{proof}

Assume first that $\cay(H,S)$ is unstable. It is sufficient to show that the S-ring 
$\A=V(G,A_{1_G})$ satisfies all the conditions in Theorem~\ref{main1}(2). i.e.,
\begin{equation}\label{eq:A}
\underline{H},\, \underline{Sa} \in \A~\text{and}~\underline{\{a\}} \notin \A. 
\end{equation}

It is clear that $\underline{Sa} \in \A$.  It has been shown above that 
$H$ and $Ha$ form an $A$-invariant partition. This implies that 
$\underline{H} \in \A$. 
Finally, due to the claim, $a_r \alpha \ne \alpha a_r$ for some 
$\alpha \in A_1$. Using also the fact that $A=G_r A_1=
A_1 G_r$, this yields that $a_r \alpha=\alpha a'_r$ holds for some 
$a' \in G$, $a \ne a'$. Then $a^{\alpha}=1^{a_r \alpha}=1^{\alpha a'_r}=a'$, 
showing that $\underline{\{a\}} \notin \A$. 
\medskip

Now assume that there is a Schurian S-ring $\A$ over $G$ satisfying all conditions in 
\eqref{eq:A}. Then $\A=V(G,B_{1_G})$ for some permutation group $B \in \Sup(G)$. 
Observe that, as $\underline{Sa} \in \A$, $B \le A$.  
Assume to the contrary that $\cay(H,S)$ is stable. Then $\alpha a_r=a_r \alpha$ for every $\alpha \in A_1$ due to the claim above, hence 
$$
\orb_{B_1}(a)=\big\{\, 1^{a_rx} : x \in  B_1\, \big\}=
\big\{\, 1^{xa_r} : x \in  B_1\, \big\}=\{a\}.
$$ 
This, however, contradicts the condition that $\underline{\{a\}} \notin \A$. 
This completes the proof of Theorem~\ref{main1}. 
%--------------------------------------------------------------------------------------------------------------%
\section{Proof of Theorems~\ref{main2} and~\ref{W-M}}\label{sec:main2+W-M}

Theorem~\ref{main2} will follow from a result of Somlai and Muzychuk~\cite{SM} 
on S-rings over abelian groups with a Sylow subgroup of prime order. We need to use 
only the case when the latter prime is $2$. 

\begin{thm}[{\cite[Propositions~3.4~and~3.5]{SM}}]\label{SM1}
Let $G$ be an abelian group of twice odd order. 
Suppose that $\A$ is an S-ring over $G$, let $Q$ be 
the least $\A$-subgroup of even order and $K$ be the largest $\A$-subgroup of odd order. Then $\A_{KQ}=\A_K \star \A_Q$.
\end{thm}

\begin{proof}[Proof of Theorem~\ref{main2}]
According to Theorem~\ref{main2}, we are given an S-ring $\A$ over $G$ such that  
$\underline{H} \in \A$ and $\underline{\{a\}} \notin \A$. 
The first condition says that $H$ is the largest $\A$-subgroup of odd order. 
Let $Q$ be the least $\A$-subgroup of even order. Then the second condition says that 
$Q > \sg{a}$, in particular, $L:=Q \cap H \ne 1$. 
Clearly, $HQ=G$, hence by Theorem~\ref{SM1}, 
$$
\A=\A_H \star \A_Q.
$$ 
By Lemma~\ref{star-is-gwp}, $\A$ is the $H/L$-wreath product.
\end{proof}

\begin{proof}[Proof of Theorem~\ref{W-M}]
Let  $\cay(H,S)$ be a connected and unstable Cayley graph. 
By Theorem~\ref{main1}, there is an S-ring  
$\A$ over the group $G$ such that $\underline{H} \in \A$, 
$\underline{\{a\}} \notin \A$ and $\underline{Sa} \in \A$. 
Since $Sa$ is a union of some basic sets of $\A$, all of which are outside $H$, 
Theorem~\ref{main2} shows that $L \le \rad(Sa)$. Consequently, as vertices of 
$\cay(H,S)$, any two $x, y \in L$ are twins.
\end{proof}

We conclude the section with a lemma, which will be used a couple of times in the next section. In the proof we shall need two further results on S-rings.

\begin{prop}[{\cite[Corollary~3.2]{SM}}]\label{SM2}
With the notations in Theorem~\ref{SM1}, suppose that $KQ \ne G$. Then 
$\A$ is the $KQ/Q$-wreath product.
\end{prop}

\begin{thm}[{\cite[Theorem~25.4]{Wbook}}]\label{W}
Let $G$ be an abelian group of composite order with a cyclic Sylow subgroup and 
$\A$ be a primitive S-ring over $G$. Then $\rk(\A)=2$.
\end{thm}

\begin{lem}\label{twice-odd}
With the notations in Theorem~\ref{SM1}, let $L=Q \cap K$. 
Then $Q \setminus L$ is a basic set. Furthermore, $\A$ is the $K/L$-wreath product. 
\end{lem}
\begin{proof}
Let $a$ be the unique involution of $G$ and $T$ be the basic set of $\A$ containing 
$a$. We show below that 
$T=Q \setminus L$. Clearly, $\underline{L} \in \A$. 
Consider the S-ring $\A_{Q/L}$. We claim that it is primitive. 
If not, then there was an $\A$-subgroup $N$ such that $L <  N < Q$. 
Since $N < Q$, $N$ cannot contain $a$ by the minimality of $Q$. 
Thus $N \le K$, so $N \le Q \cap K=L$, contradicting the assumption that 
$N > L$. By Theorem~\ref{W}, 
$\rk(\A_{Q/L})=2$. This then combined with \eqref{eq:cap} yield the 
existence of a positive number $\ell$ such that 
$$
|Lx \cap T|=\ell~\text{for every}~x \in Q \setminus L. 
$$
 
On the other hand, $\A_{KQ}=\A_K \star \A_Q$ by Theorem~\ref{SM1}. 
This shows that $La \subset T$, so $\ell=|L|$, i.e., $T=Q \setminus L$.   

Let $X \in \S(\A)$ be an arbitrary basic set outside $K$. We have to show 
that $L \le \rad(X)$.  If $X \not\subseteq KQ$, then due to Proposition~\ref{SM2}, 
$\rad(X) \ge Q > L$. If $X \subseteq KQ$, then $L \le \rad(X)$ follows from 
Lemma~\ref{star-is-gwp}.
\end{proof}
%--------------------------------------------------------------------------------------------------------------%
\section{Proof of Theorem~\ref{main3}}\label{sec:main3}

For this section we set the following assumptions. 

\begin{hypo}\label{hypos}
$H$ is an abelian group with a unique involution $b$ and 
$H_0 < H$ is the unique subgroup of $H$ of order $|H|/2$ and $|H| >2$. 
Furthermore, 
\begin{quote}
$\A$ is an S-ring over $G=H \times \sg{a}$, where $\sg{a} \cong \Z_2$ 
such that $\underline{H} \in \A$. \\ [+1ex] 
$T$ is the basic set of $\A$ containing $a$. \\ [+1ex]
$K$ is the largest $\A$-subgroup of odd order.
\end{quote}
\end{hypo}

The proof of Theorem~\ref{main3} will be given in the end of the section 
following four preparatory lemmas. 

The following simple fact will be used often. 
If $A, B \le G$ are any subgroups and $S \subset G$ is any subset, then 
\begin{equation}\label{eq:rad}
AB/B \le \rad(S/B) \implies A \le \rad(SB). 
\end{equation}

\begin{lem}\label{L1}
Assuming Hypothesis~\ref{hypos}, suppose that $\underline{Kab} \in \A$ and $La \subseteq T$ for some $\A$-subgroup $L$, $L \le H_0$. Then
$$ 
L \le \bigcap_{X \in \mathcal{S}(\A) \atop X \not\subseteq H \cup Kab}\rad(X). 
$$
\end{lem}
\begin{proof}
Fix a basic set $X  \in \mathcal{S}$ such that $X \not\subseteq H \cup Kab$. 
We show below that $L \le \rad(X)$. As $\underline{L} \in \A$, there is a positive number 
$\ell$ such that $|X \cap Lx|=0$ or $\ell$ for every $x \in G$, see \eqref{eq:cap}.
As $X \not\subseteq H$, $X$ can be expressed as 
$$
X=X_1a \cup X_2ab \cup X_3a \cup X_3ab,
$$  
where $X_1, X_2$ and $X_3$ are pairwise disjoint subsets of $H_0$. 

Assume first that $X_1 \cup X_2=\emptyset$.
Then $b \in \rad(X)$, hence $Q \le \rad(X)$, where $Q$ is the 
least $\A$-subgroup containing $b$. Let us consider the S-ring 
$\A_{G/Q}$. Then $G/Q$ has twice odd order and $T/Q$ is a basic set of 
$\A_{G/Q}$ containing the unique involution of $G/Q$. It follows from 
Lemma~\ref{twice-odd} that  

$$
T/Q=\sg{T}Q/Q \setminus H/Q~\text{and}~
\sg{T}Q/Q \cap H/Q \le \rad(X/Q).
$$ 
The group $\sg{T}Q/Q \cap H/Q =(\sg{T} \cap H)Q/Q$. Using \eqref{eq:rad}, 
we obtain $\sg{T} \cap H \le \rad(XQ)=\rad(X)$. As $La \subseteq T$, 
$L \le \sg{T} \cap H$, so $L \le \rad(X)$. 

Now assume that $X_1 \cup X_2 \ne \emptyset$,
Then $X^{[2]}=(X_1 \cup X_2)^{(2)}$. Due to Theorem~\ref{S-multi}(2), the 
latter set is an $\A$-set, which is clearly contained in $K$. 
As $|K|$ is odd, there is an integer 
$m$ such that $\gcd(m,|K|)=1$ and $2m \equiv 1 \!\!\pmod {|K|}$. 
Using Theorem~\ref{S-multi}(1), we conclude that 
$X_1 \cup X_2=(X_1 \cup X_2)^{(2m)}$ is also an $\A$-set. 
If $X_2 \ne \emptyset$, then 
$X \cap Kab \ne \emptyset$, hence 
$X \subseteq Kab$. This is impossible by our assumption that $X \not\subseteq Kab$, 
thus $X_2=\emptyset$ and $X_1 \ne \emptyset$. Then $\underline{X} \cdot \underline{X_1^{(-1)}} \in \A$. We have $\underline{X} \cdot \underline{X_1^{(-1)}}=
\sum_{x \in G}\alpha_x x$ for some non-negative integers $\alpha_x$'s.  
It is easy to see that $\alpha_a=|X_1|$. Also, $\alpha_y=\alpha_a$ for every $y \in T$ 
because $T$ is a basic set and $a \in T$. In particular, as $La \subseteq T$, 
we obtain that 
$$
\sum_{y \in La}\alpha_{y}=|X_1| \cdot |L|.
$$

Now fix $x \in X_1$. Denote by $\nu_x$ the number of elements $x' \in X$ such that 
$x' x^{-1} \in La$. We find that $\nu_x=|X \cap Lax|=\ell$ because $ax \in X$. 
Then we can write that 
$$
|X_1| \cdot |L|=\sum_{y \in La}\alpha_{y}=\sum_{x \in X_1}\nu_x=|X_1| \cdot \ell.
$$
This shows that $\ell=|L|$, so $L \le \rad(X)$. 
\end{proof}

\begin{lem}\label{L2}
Assuming Hypothesis~\ref{hypos}, suppose that $ab \in T$ and $La \subseteq T$ 
for some $\A$-subgroup $L$, 
$L \le H_0$. Then $\A$ is the $H/L$-wreath product.
\end{lem}
\begin{proof}
Assume to the contrary that there is a basic set $X$, $X \not\subseteq H$ and 
$L \not\le \rad(X)$. Due to \eqref{eq:cap}, there is a constant $\ell$, 
$0 < \ell < |L|$ such that $|X \cap Lx|=0$ or $\ell$ for every $x \in G$. Since $|L|$ is odd, it is possible to choose $X$ so that $\ell < |L|/2$. 

As $X \not\subseteq H$, $X=X_1a \cup X_2ab \cup X_3a \cup X_3ab$, where 
$X_1, X_2$ and $X_3$ are pairwise disjoint subsets of $H_0$. 
If $X_1 \cup X_2 = \emptyset$, then the argument, used in the proof of the 
previous lemma, yields that $L \le \rad(X)$. This is impossible, hence 
$X_1 \cup X_2 \ne \emptyset$.

Consequently, $X_1 \cup X_2$ is a non-empty $\A$-set, and 
the product $\underline{X} \cdot \underline{X_1^{(-1)} \cup X_2^{(-1)}}$  
belongs to $\A$. Write it as $\sum_{x \in G}\alpha_x x$.  
It is easy to see that $\alpha_a=|X_1|$ and $\alpha_{ab}=|X_2|$. Since $ab \in T$, 
$\alpha_a=\alpha_{ab}$, so $|X_1|=|X_2|$. As $La \subseteq T$, we obtain 
$$
\sum_{y \in La}\alpha_{y}=|X_1| \cdot |L|.
$$

Now fix $x \in X_1 \cup X_2$ and denote by $\nu_x$ the number of 
elements $x' \in X$ such that $x' x^{-1} \in aL$. 
Notice that $\nu_x=|X \cap Lax|$, and so $\nu_x=0$ or $\ell$ for every 
$x \in X_1 \cup X_2$. Then we can write 
$$
|X_1| \cdot |L|=\sum_{y \in La}\alpha_{y}=\sum_{x \in X_1 \cup X_2}\nu_x   
\le (|X_1|+|X_2|) \cdot \ell=|X_1| \cdot 2\ell.
$$
This contradicts our assumption that $\ell < |L|/2$.
\end{proof}

\begin{lem}\label{L3}
Assuming Hypothesis~\ref{hypos}, suppose that $T=La \cup Lab$ for some subgroup $L \le H_0$, $L \ne 1$, and 
$L$ contains no non-trivial $\A$-subgroup. Then $\A$ is the $H/M$-wreath product, 
where $M=\sg{b,L}$.
\end{lem}
\begin{proof}
Observe that $\sg{T}=\sg{a,b} L$. 
Thus $\sg{T} \cap H=M$, in particular, $\underline{M} \in \A$ because both 
$\sg{T}$ and $H$ are $\A$-subgroups. 
Let $N$ be a minimal non-trivial $\A$-subgroup contained 
in $M=\sg{L,b}$. Then $\A_N$ is a primitive S-ring.  As $N  \not\le L$, 
$\sg{b}$ is a Sylow $2$-subgroup of $N$. By Theorem~\ref{W}, $\A_N$ has rank 
$2$, and we have that $N^\#$ is a basic set.

Consider the S-ring $\A_{G/N}$.  Then $G/N$ has twice odd order and 
$T/N$ is the basic set containing the unique involution. 
It follows from Lemma~\ref{twice-odd} that
$$
T/N=\sg{T}/N \setminus H/N~\text{and}~\sg{T}/N \cap H/N \le \rad(X/N),
$$
where $X  \in \S(\A)$, $X \not\subseteq H$. 
The group $\sg{T}/N \cap H/N=M/N$, and by \eqref{eq:rad}, 
$M \le \rad(NX)$. This shows that it is sufficient to show that 
$N \le \rad(X)$ for every basic set $X \in \S(\A)$, $X \not\subseteq H$.  

Assume to the contrary that there is a basic set $X$ such that $X \not\subseteq H$ 
and $N \not\le \rad(X)$. Due to Eq.~\eqref{eq:cap}, there is a constant $\ell$, 
$0 < \ell < |N|$ such that $|X \cap Nx|=0$ or $\ell$ for every $x \in G$. 
It is possible to choose $X$ such that $\ell \le |N|/2$.

As $X \not\subset H$, $X=X_1a \cup X_2ab \cup X_3a \cup X_3ab$, where 
$X_1, X_2$ and $X_3$ are pairwise disjoint subsets of $H_0$. 
Let us consider the product $\underline{X} \cdot \underline{X^{(-1)}}$, which is 
in $\A$. Write it as $\sum_{x \in G}\alpha_x x$.  
It is easy to see that $\alpha_{b}=2|X_3|$. As $N^\#$ is a basic set, we obtain 
$$ 
\sum_{y \in N^\#}\alpha_{y}=2|X_3| \cdot (|N|-1).
$$
Now fix $x \in X$. If $\nu_x$ denotes the number of elements $x' \in X$ such that 
$x' x^{-1} \in N^\#$, then we find that $\nu_x=|X \cap Nx|-1=\ell-1$, and so we obtain that 
\begin{equation}\label{eq:1}
2|X_3| \cdot (|N|-1)=\sum_{y \in N^\#}\alpha_{y}=
\sum_{x \in X}\nu_x=|X| \cdot (\ell-1).
\end{equation}

This combined with the fact that $|X|=2|X_3|+|X_1|+|X_2|$ and the
assumption that $\ell \le |N|/2$ yield $2|X_3| < |X_1|+|X_2|$.
In particular, $X_1 \cup X_2 \ne \emptyset$. 
Let $\underline{X} \cdot \underline{(X_1^{(-1)}+X_2^{(-1)})}=\sum_{x \in G}\beta_x x$. Computing the value $\sum_{y \in Na}\beta_y$ in two ways as in  
the proof of Lemma~\ref{L2}, one can deduce that 
$$
|X_1| \cdot |N|=|X_1| \cdot 2 \ell=|X_1| \cdot 2.
$$
Thus $N=\sg{b}$. This forces that $\ell=1$, and thus $X_3=\emptyset$ follows from 
\eqref{eq:1}.  
We have shown above that $M \le \rad(NX)$, and so for $x \in X$, 
$$
|M|=|Mx \cap  NX|=|Mx \cap (X \cup Xb)|=
|Mx \cap X|+|Mx \cap Xb|=2|Mx \cap X|.
$$
The third equality is true because $X \cap Xb=\emptyset$ and 
the fourth equality follows as $b \in M$. This shows that there are 
two elements $y, z \in  Mx \cap X$. Recall that $X=X_1a \cup X_2ab$. 
If $y, z \in X_1a$ or $X_2ab$, then $y z^{-1} \in M \cap (X_1 \cup X_2)$.  
If $y \in X_1a$ and $z \in X_2ab$, then $yb \in Mx \cap X_1ab$, and 
$1 \ne yb z^{-1} \in M \cap \sg{X_1  \cup X_2}$ ($yb \ne z$ because 
$X_1 \cap X_2=\emptyset$).  
Noticing that $X_1 \cup X_2 \subseteq K$, we find that $K \cap M \ne 1$, 
in particular, $K \cap M$ is a non-trivial $\A$-subgroup contained in $L$.  
This contradicts our initial assumption that no such subgroup exists. 

\end{proof}

In our last lemma before the proof of Theorem~\ref{main3} we describe 
the basic set $T$ when $H$ is a cyclic group.

\begin{lem}\label{L4}
Assuming Hypothesis~\ref{hypos}, suppose that $H$ is a cyclic group. Then 
\begin{equation}\label{eq:T}
T \in \big\{ La,~La \cup Lab,~Ma \cup (M \setminus L)ab : 1 \le L < M \le H_0 \big\}. 
\end{equation}
\end{lem}
\begin{proof}
We proceed by induction on $|H_0|$. 
Suppose first that $|H_0|=p$ for a prime $p$.  For every integer 
$k$ such that $\gcd(k,2p)=1$, $a^k=a$, and thus $T^{(k)}=T$ due to 
Theorem~\ref{S-multi}(1). It follows that $T$ is one of the following sets:
$$ 
\{a\},~\{a,ab\},~\{a\} \cup H_0^\#ab,~\{a\} \cup H_0ab,~
H_0a,~H_0a \cup \{ab\},~H_0a \cup H_0^\#ab,~H_0a \cup H_0ab.
$$   
Thus \eqref{eq:T} holds unless $T=\{a\} \cup H_0^\#ab$ or 
$\{a\} \cup H_0ab$ or $H_0a \cup \{ab\}$. In each of the latter cases, 
$H_0=\sg{T^{[2]}}$, so $\underline{H_0} \in \A$ by Theorem~\ref{S-multi}(2). 
Then, however, $|T \cap H_0a| \ne |T \cap H_0ab|$, contradicting the identity in  ~\eqref{eq:cap}. 
This shows that the lemma holds if $|H_0|$ is a prime. 

Now assume that $|H_0|$ is a composite number.
Let $R=\rad(T)$. If $R \ne 1$ and $|R|$ is odd, then the lemma follows from  the induction hypothesis applied to $\A_{G/R}$. Whereas if $|R|$ is even, then the lemma follows from Lemma~\ref{twice-odd} applied to $\A_{G/R}$. 
For the rest of the proof let $R=1$. We are going to show that 
$T=M\{a,ab\} \setminus \{ab\}$ for some $1 \le M \le H_0$, in particular,  
\eqref{eq:T} holds in this case as well.  

Fix a prime divisor $p$ of $|H_0|$ and consider the set $T^{[p]}$. 
Since $R=1$ and $H_0$ is a cyclic group, it follows that $T^{[p]} \ne \emptyset$. 
Let $N=\sg{T^{[p]}}$. It is clear that $N < G$ and $\underline{N} \in \A$ by 
Theorem~\ref{S-multi}(2). 
If $\sg{a,b} \le N$, then the 
induction hypothesis can be applied to $\A_N$, and this yields 
$T=M\{a,ab\} \setminus \{ab\}$ for some $1 \le M \le H_0$.
Therefore, we may assume that 
$|\sg{a,b} \cap N|=2$.
 
Now if $a \in N$, then we can apply Lemma~\ref{twice-odd} to $\A_N$ and conclude 
that $T=\{a\}$  because $\rad(T)=1$. 

Therefore, $ab \in N$ but $a \notin N$. We show below that these conditions give rise to a 
contradiction. Using Lemma~\ref{twice-odd} and the fact that 
$\underline{H} \in \A$, we find that the basic set of $\A$ 
containing $ab$ is equal to $Lab$ for some subgroup $L \le H_0$. 
Write $T$ as 
$$
T=T_1a \cup T_2ab \cup T_3a \cup T_3ab,
$$ 
where $T_1, T_2$ and $T_3$ are pairwise disjoint subsets of $H_0$. Recall that $K$ is 
the largest $\A$-subgroup of odd order.  As $\underline{L} \in \A$, $L \le K$, and 
hence $KT=KLab=Kab$ is an $\A$-subset. It is clear that $Kab \cap T=\emptyset$. 
Using also that $T_1 \cup T_2 \subseteq K$, we find that $T_2=\emptyset$. 
On the other hand, the condition that $a \notin N$ shows that $P \le \rad(T_1)$ and 
$P \le \rad(T_3)$, where $P$ is the subgroup of $H_0$ of order $p$. 
We conclude that $P \le \rad(T)=R$, contradicting our assumption that $R=1$. 
\end{proof}

We are ready to prove Theorem~\ref{main3}.

\begin{proof}[Proof of Theorem~\ref{main3}]
Let us keep all the notations $H, H_0, G, a, b, \A, T$ set in Hypothesis~\ref{hypos}, and 
assume, in addition, that $H$ is a cyclic group and $\underline{\{a\}} \notin \A$. 
Define the subgroup 
$$
V=\bigcap_{X \in \mathcal{S}(\A) \atop X \cap H_0a \ne \emptyset}\rad(X \cap H_0a).
$$
We have to show that $V \ne 1$ provided that $\{a,ab\}$ is not a basic set of $\A$.  
We distinguish three cases according to the possibilities for $T$ mentioned in  Lemma~\ref{L4}.
\medskip

\noindent{\bf Case~1.} $T=La$, $1< L \le H_0$. 
\medskip

Let $X \in \S(\A)$ such that $X \cap H_0a \ne \emptyset$. 
It is sufficient to show that $L \le \rad(X)$. Notice that $\sg{T,K}=\sg{a,K}$ is an $\A$-subgroup.  Applying Lemma~\ref{twice-odd} to $\A_{\sg{a,K}}$, we obtain that 
$L \le \rad(X)$ if $X \subseteq Ka$. 
Assume that $X \not\subseteq Ka$. 
Let $T'$ be the basic set containing $ab$. 
By Lemma~\ref{L4}, $T'=Mab$ or $Nab \cup (N \setminus M)a$ 
for some subgroups $1 \le M < N \le H_0$. 
In the former case $M \le K$ and $\underline{Kab} \in \A$. As 
$X \not\subseteq H \cup Kab$, $L \le \rad(X)$ due to Lemma~\ref{L1}. 
In the latter case $\sg{T'} \setminus (H \cup T')=Ma$, implying that 
$\underline{Ma} \in \A$. Thus $L \le M$, and so 
$Lab \subseteq T'$. As $X \not\subseteq H \cup Ka$, $L \le \rad(X)$ due to Lemma~\ref{L1} applied for the basic $T'$ and the element $ab$. 
%If $X \not\subseteq Ka$ and $X \cap H_0a \ne \emptyset$, then 
%$b \in \rad(X)$. Thus $Q \le \rad(X)$, where $Q$ is the least $\A$-subgroup containing %$b$. Then $LQ/Q \le \rad(X/Q)$ follows from Lemma~\ref{twice-odd} when applied to $%\A_{G/Q}$.  This yields that $L \le \rad(QX)=\rad(X)$. 
\medskip

\noindent{\bf Case~2.} $T=La \cup Lab$, $1 \le L \le H_0$.  
\medskip

If $L=1$, then $\{a,ab\}$ is a basic set. Assume that $L \ne 1$. 
Then it follows from Lemmas~\ref{L2} and \ref{L3} that 
$\A$ is the $H/N$-wreath product, where $N=L$ or $N=\sg{L,b}$.  
Thus $N \cap H_0=L$, implying that $1 < L \le V$, in particular, $V \ne 1$. 
\medskip

\noindent{\bf Case~3.} $T=La \cup (L \setminus M)ab$,  $1 \le M < L \le H_0$.
\medskip

Then $\sg{T} \setminus (H \cup T)=Mab$, implying that $\underline{Mab} \in \A$.  
By Lemma~\ref{L4}, the basic set containing $ab$ is equal to the 
coset $Nab$ for some subgroup $1 \le N \le M$. 

If $N=1$, then $\{ab\}$ is a basic set. Then so is $Tab=Lb \cup (L \setminus M)$. 
Observe that if $X \in \S(\A)$ such that $X \cap H_0a \ne \emptyset$, then 
$Xab$ is also basic, it is contained in $H$ and has non-empty intersection with 
$H_0b$. By Lemma~\ref{twice-odd}, $L \le \rad(Xab)$, so $L \le \rad(X)$, and 
this yields that $1 < L \le V$.

Now assume that $N \ne 1$. Let $X \in \mathcal{S}(\A), X \not\subseteq H$. If $X \subseteq Kab$, then applying Lemma~\ref{twice-odd} to $\A_{\sg{Nab,K}}$, 
we obtain that $N \le \rad(X)$. If $X \not\subseteq Kab$, then $N \le \rad(X)$ holds by 
Lemma~\ref{L1}.  All these yield that $1 < N \le V$. 
\end{proof}
%--------------------------------------------------------------------------------------------------------------%
\section{Isomorphisms of circulant graphs}\label{sec:iso}

Muzychuk~\cite{M} derived an effective criterion for the isomorphisms between 
two circulant graphs in terms of generalised multipliers. 
In this section we give a short overview of his result. For our purpose it will 
be enough to consider the particular case when the order of the two circulant graphs is 
equal to $2p^e$ for some odd prime $p$ and integer $e \ge 1$. 

\begin{defi}\label{key}
Let $p$ be a prime and $e \ge 1$ be an integer.
The \emph{key space} $\K_{p^e}$ consists of the $e$-tuples 
$\k=(k_1,\ldots,k_e)$ of integers such that 
\smallskip

\begin{enumerate}[(1)]
\setlength{\itemsep}{0.25\baselineskip}
\item $\forall 1 \le i \le e :~0 \le k_i \le i-1$,
\item $\forall 2 \le i \le e :~k_{i-1} \le k_i$.
\end{enumerate}
\smallskip
\end{defi}

The elements in $\K_{p^e}$ are called \emph{primary keys}.

\begin{defi}\label{key-part}
Let $H  \cong \Z_{2p^e}$ for an odd prime $p$ and $e \ge 1$ and   
$\k=(k_1,\ldots,k_e) \in \K_{p^e}$. Then the \emph{key partition} $\Pi(\k)$ is 
the partition of $H$ defined as 
$$
\Pi(\k)=\big\{\, \{1_G\}, P_{k_i}x : x \in G^\#, 
o(x)=p^i~\text{or}~2p^i, 1 \le i \le e\, \big\},
$$
where for every $i$, $1 \le i \le e$, 
$P_{k_i} < H$ is the subgroup of order $p^{k_i}$. 
\end{defi}

Note that, for every primary key $\k=(k_1,\ldots,k_e) \in \K_{p^e}$, 
the $(e+1)$-tuple $(p^{k_1},\ldots,p^{k_e};\Pi_{=})$ is an S-system, 
where $\Pi_{=}$ is the partition of $[e]$ consisting of the singleton subsets $\{i\}$, 
$1 \le i \le e$. 
Let $H \cong \Z_{2p^e}$. 
Denoting by $H_0$ be the unique subgroup of $H$ of order 
$p^e$, the following equality holds (see also~\cite[part (1) in Proposition~4.11]{M}): 
\begin{equation}\label{eq:class=basic-set}
\big\{ X \in \Pi(\k) : X \subseteq H_0 \big\}=\S(\A),~\text{where}~ 
\A=\S_{H_0}(p^{k_1},\ldots,p^{k_e};\Pi_{=}).   
\end{equation}

There is a natural partial order $\le$ on the key space $\K_{p^e}$ defined by 
$$
\k \le \k'~\stackrel{\text{def.}}{\iff}~\forall 1\le  i \le e:~k_i \le k'_i, 
$$
where $\k, \k' \in \K_{p^e}$, $\k=(k_1,\ldots,k_e)$, and $\k'=(k'_1,\ldots,k'_e)$. 
The poset $(\K_{p^e},\le)$ is a lattice, whose \emph{meet} $\wedge$ and 
\emph{join} $\vee$ operators are defined as  
$$
\forall 1 \le i \le t:~
(\k \wedge \k')_i=\min(k_i,k'_i)~\text{and}~(\k \vee \k')_i=\max(k_i,k'_i).
$$

A non-empty subset $S \subseteq H, H \cong \Z_{2p^e}$ is said to be a 
\emph{$\Pi(\k)$-subset} if $S$ can be expressed as the union of some classes of 
$\Pi(\k)$. If $S$ is a $\Pi(\k')$-subset for 
another primary key $\k'$, then it is also a $\Pi(\k \vee \k')$-subset. 
This observation allows us to consider the greatest element $\k$ 
of the poset $(\K_{p^e},\le)$ for which $S$ is a $\Pi(\k)$-subset.

\begin{defi}\label{key-S}
Let $H  \cong \Z_{2p^e}$ for an odd prime $p$ and $e \ge 1$.  
The \emph{key} of a non-empty subset $S \subseteq H$ is 
the greatest primary key $\k \in \K_{p^e}$ in the poset $(\K_{p^e},\le)$ 
for which $S$ is a $\Pi(\k)$-subset.
\end{defi}

\begin{defi}\label{gen-multi}
Let $p$ be a prime and $e \ge 1$ be an integer.
A \emph{generalised multiplier} of $\Z_{p^e}$ is a 
an $e$-tuple $\vec{m}=(m_1,\ldots,m_e)$ of positive integers 
such that $\gcd(m_i,p)=1$ for every $1 \le i \le e$.
\end{defi}

The set of all generalised multipliers of $\Z_{p^e}$ is denoted by 
$\Z_{p^e}^{**}$. 

\begin{defi}\label{f_m}
Let $p$ be a prime and $e \ge 1$ be an integer, and 
let $\vec{m}=(m_1,\ldots,m_t) \in \Z_{p^e}^{**}$. Define $f_{\vec{m}}$ to be 
the permutation of $\Z_{p^e}$ acting by the rule
$$
\forall x\in \Z_{p^e} :~x^{f_{\vec{m}}}=\sum_{i=0}^{e-1} m_{e-i}x_ip^i,
$$
where $x=\sum_{i=0}^{e-1}x_i p^{i}$ is the \emph{$p$-adic expansion} of $x$, 
i.e., $0 \le  x_i \le p-1$ for every index $i$, $0 \le i \le e-1$.
\end{defi}

\begin{defi}\label{gen-multi-k}
Let $p$ be a prime and $e \ge 1$ be an integer.
For a primary key $\k \in \K_{p^e}$, define $\Z_{p^e}^{**}(\k)$ to be the set of 
all generalised multipliers $\vec{m}=(m_1,\ldots,m_e)$ for which 
$$
\forall 2 \le i \le e:~m_i \equiv m_{i-1}\!\!\!\!\pmod {p^{i-1-k_i}}.
$$
\end{defi}

\begin{prop}[{\cite[case~(2) in Proposition~2.4]{M}}]\label{M1}
Let $H=\sg{h}  \cong \Z_{2p^e}$ for an odd prime $p$ and $e \ge 1$. 
Let $\k \in \K_{p^e}$ and $\vec{m}=(m_1,\ldots,m_e) \in \Z_{p^e}^{**}(\k)$ 
be a generalized multiplier such that $m_i$ is odd for every $i$, $1 \le i \le e$. 
Let $\varphi_{\vec{m}}$ be the permutation of $H$ defined as 
\begin{equation}\label{eq:f}
\forall x \in \Z_{p^e} : (h^{2x})^{\varphi_{\vec{m}}}=
h^{2x^{f_{\vec{m}}}}~\text{and}~
(h^{2x+1})^{\varphi_{\vec{m}}}=h^{2x^{f_{\vec{m}}}+1}.
\end{equation}
Then for every class $X \in \Pi(\k)$, if $|\sg{X}|_p=p^i$, then 
$X^{\varphi_{\vec{m}}}=X^{(m_i)}$.
\end{prop}

We are ready to formulate Muzychuk's criterion. 

\begin{thm}[{\cite[Theorem~1.1]{M}}]\label{M2}
Let $H=\sg{h}  \cong \Z_{2p^e}$ for an odd prime $p$ and $e \ge 1$. Then 
two graphs $\cay(H,S)$ and $\cay(H,S')$ are isomorphic 
if and only if 
\begin{enumerate}[(1)]
\item $S$ and $S'$ have the same key, say $\k$.
\item There exists a generalised multiplier $\vec{m} \in \Z_{p^e}^{**}(\k)$ such that 
$S'=S^{\varphi_{\vec{m}}}$, where $\varphi_{\vec{m}}$ is the permutation of $H$ defined in \eqref{eq:f}.
\end{enumerate}
\end{thm}
%--------------------------------------------------------------------------------------------------------------%
\section{Proof of Theorem~\ref{main4}}\label{sec:main4}

We need one more lemma about S-rings.

\begin{lem}\label{calB}
Let $H=E \times F$ be an abelian group such that 
$E=\sg{u,v} \cong \Z_2^2$ and $|F|$ is odd, and suppose that 
$\A$ is a Schurian S-ring over $H$ such that 
$\underline{F}, \underline{\sg{F, v}}, \in \A$ and 
$\underline{\{u,uv\}} \in \S(\A)$. 
Let $X \in \S(\A)$, $X \not\subset \sg{F, v}$. Then 
$$
|X \cap Fu|=|X \cap Fuv|.
$$
Furthermore, both sets $X \cap Fu$ and $X \cap Fuv$ are the basic sets of 
a Schurian S-ring $\B$ over $H$, for which $\A \subset \B$. 
\end{lem}
\begin{proof}
As $\A$ is Schurian, $\A=V(H,A_{1_H})$ for a group $A \in \Sup(H)$.  
Let $K$ be the kernel of the action of $A$ on the set $[H:F]$ consisting of the 
$F$-cosets in $H$, and let $B=K\sg{u_r, v_r}$. 
Clearly, $B \in \Sup(H)$. Let $\B=V(H,B_{1_H})$. 

Since $B < A$, it follows that $\B \supset \A$. 
The sets $\{u,uv\}/F$ and $X/F$ are basic sets of the S-ring $\A_{H/F}$ whose 
intersection is non-empty. Thus, they are equal, implying that there are 
elements $x_1 \in X \cap Fu$ and $x_2 \in X \cap Fuv$.   

Let $X_i$ be the basic set of $\B$ containing $x_i$ for $i=1,2$. 
Clearly, $X_1$ and $X_2$ are contained in $X$ and belong to different $F$-cosets. 
Observe that $A_{1_H} \cap A_{x_1} \le K$. It follows from this that  
$A_{1_H} \cap A_{x_1}=B_{1_H} \cap B_{x_1}$. 
Also, $|B|=|K| |E|=|A|/2$, and therefore, $|B_{1_H}|=|A_{1_H}|/2$. 
These together with the orbit-stabilizer lemma yield 
$$
|X_1|=\frac{|B_{1_H}|}{|B_{1_H} \cap B_{x_1}|}=
\frac{|A_{1_H}|}{2 |A_{1_H} \cap A_{x_1}|}=|X|/2.
$$
The same argument shows that $|X_2|=|X|/2$, and so $|X_1|=|X_2|$ and 
$X_1=X \cap Fu$ and $X_2=X \cap Fuv$. 
\end{proof}

We are ready to prove Theorem~\ref{main4}. 

\begin{proof}[Proof of Theorem~\ref{main4}]
Let $H=\sg{h} \cong \Z_{2p^e}$ for some odd prime $p$ and $e \ge 1$ and  
$\cay(H,S)$ be the graph described in Theorem~\ref{main4}, i.e., 
$\cay(H,S)$ is a connected and non-bipartite. 
Let $b$ be the unique involution of $H$ and $H_0$ the unique subgroup of 
$H$ of order $p^e$. 
 
In view of Remark~\ref{to:main4}, we may assume that $\cay(H,S)$ is 
unstable. To settle Theorem~\ref{main4}, we have to show that, 
\begin{equation}\label{eq:cond1}
(S \cap H_0)h=S \cap H_0~\text{for some}~h \in H_0, h \ne 1_H,~\text{or}~
\end{equation}
\begin{equation}\label{eq:cond2}
\cay(H,S) \cong \cay(H,Sb).
\end{equation}

Let $G=H \times \sg{a}$, where $\sg{a} \cong \Z_2$, and let 
$\A$ be the intersection of all Schurian S-rings over $G$ such that 
$\underline{H}, \underline{Sa} \in \A$. Then $\A$ is Schurian and it follows from Theorem~\ref{main1} that $\underline{\{a\}} \notin \A$. 
Let 
$$
V=\bigcap_{X \in \mathcal{S}(\A) \atop X \cap H_0a \ne \emptyset}\rad(X \cap H_0a).
$$
Then $S \cap H_0 \ne \emptyset$ because $\cay(H,Sa)$ is non-bipartite. 
Thus $Sa \cap H_0a \ne \emptyset$, hence   
$(Sa \cap H_0a)h=Sa \cap H_0a$ for evey $h \in V$. 
This shows that \eqref{eq:cond1} holds if $V \ne 1$. 
For the rest of the section we assume that $V=1$. 
\medskip

Let $\k=(k_1,\ldots,k_e)$ be the key of $S$. 
In view of Theorem~\ref{M2}, we derive \eqref{eq:cond2} by 
finding a generalized multiplier $\vec{m} \in \Z_{p^e}^*(\k)$ for which  
$S^{\varphi_{\vec{m}}}=Sb$. 

Observe that $\{a,ab\}$ is a basic set of $\A$ due to Theorem~\ref{main3}. 
Thus $E:=\sg{a,b}$ is an $\A$-subgroup, and by Theorem~\ref{P2}, 
$$
\A_{G/E} \cong \S_{G/E}(d_1,\ldots,d_e;\Pi)~\text{for some 
S-system}~(d_1,\ldots,d_e;\Pi).
$$ 

\begin{claim1} $\forall 1 \le i \le e: p^{k_i}=(d_i)_p$. 
\end{claim1}
\begin{proof}
For $1 \le i \le e$, put $p^{k'_i}=(d_i)_p$ and 
let $\k'=(k'_1,\ldots,k'_e)$. Then $\k' \in \K_{p^e}$ and 
as $\underline{Sa} \in \A$, $S$ is a $\Pi(\k')$-subset.  
Thus $\k' \le \k$ due to Definition~\ref{key-S}.

In view of \eqref{eq:class=basic-set}, we can define the S-ring $\A'$ over $G$ as 
\begin{equation}\label{eq:A'}
\A'=\S_{H_0}(k_1,\ldots,k_e;\Pi_=) \otimes \Z E.
\end{equation} 
It can be easily seen that $\underline{H}, \underline{Sa} \in \A'$. 
Also, $\A'$ is Schurian because it is the tensor product of 
two Schurian S-rings (see Theorem~\ref{P1}), and so 
$\A \subseteq \A'$. Then $\A_{G/E} \subseteq \A'_{G/E}$, hence every basic set of $\A_{G/E}$ is union of some basic sets of $\A'_{G/E}$. Fix a number $i$, $1 \le i \le e$. 
Let $X \in \S(\A_{G/E})$ be a basic set such 
that $X$ contains an element of order $p^i$ and let 
$X' \in \S(\A'_{G/E})$ such that $X' \subseteq X$. 
If $d_i \ne p^{i-1}(p-1)$, then $|X|=d_i$, see \eqref{eq:Delta}.
Using also that $|X'|=p^{k_i}$, we find that 
$p^{k_i} \le (d_i)_p=p^{k'_i}$. Since $k_i \le i-1$, see Definition~\ref{key}(1), 
$p^{k_i} \le p^{k'_i}$ also holds when $d_i=(p-1)p^{i-1}$, and we conclude that  
$\k \le \k'$. As $\k \ge \k'$ also holds, we have $\k=\k'$. 
\end{proof}

For $1 \le i \le e$, let $P_i$ be the subgroup of $H$ of order $p^i$. 

\begin{claim2}
Let $X \in \S(\A)$ be a basic set such that $X \not\subseteq H$ and $X \ne \{a,ab\}$, 
and let $x \in X$ be an element such that $o(x)=2p^i$.  
Then $P_{k_i}x \subseteq X$. 
\end{claim2}
\begin{proof}
The set $P_{k_i}x$ is a basic set of the S-ring $\A'$ defined in \eqref{eq:A'}. 
As $\A \subseteq \A'$, the claim follows.  
\end{proof} 

For positive integers $i$ and $j$ such that $\gcd(i,j)=1$, 
let $o_i(j)$ denote the \emph{multiplicative order} of $j$ mod $i$, i.e., it is the 
smallest positive integer $k$ such that $j^k \equiv 1\!\!\pmod i$. 
In the proof of our next claim we also need the following elementary fact from 
number theory: Let $i, m$ be positive integers such 
that $\gcd(p,m)=1$. Then  
\begin{equation}\label{eq:mod}
\gcd(p,o_{p^{i+1}}(m))=1~\implies~o_{p^i}(m)=o_{p^{i+1}}(m). 
\end{equation}

\begin{claim3} 
There exists a generalised multiplier 
$\vec{m}=(m_1,\ldots,m_e) \in \Z_{p^e}^{**}(\k)$ such that  
$$
\forall 1 \le i \le e:~\gcd(m_i,2p)=1~\text{and}~
o_{p^i}(m_i)=(d_i)_{p'}.
$$
\end{claim3}
\begin{proof}
We define the entries of $\vec{m}$ recursively. 
Choose first $m_e$ to be a positive integer such that 
$\gcd(m_e,2p)=1$ and $o_{p^e}(m_e)=(d_e)_{p'}$. 
The existence of $m_e$ is guaranteed by the Chinese remainder theorem. 
Suppose that $m_i$ is already defined for an integer $i$, $2 \le i \le e$. 
If $k_i=i-1$, then define $m_{i-1}$ to be any positive integer such that  
$\gcd(m_{i-1},2p)=1$ and $o_{p^{i-1}}(m_{i-1})=(d_{i-1})_{p'}$, 
and if $k_i < i-1$, then let $m_{i-1}=m_i$.

It follows at once that $\vec{m} \in \Z_{p^e}^{**}$ (see Definition~\ref{gen-multi}). 
Moreover, $m_i \equiv m_{i-1}\!\!\pmod  {p^{i-1-k_i}}$ whenever $2 \le i \le e$, 
which means that $\vec{m} \in \Z_{p^t}^{**}(\k)$ (see 
Definition~\ref{gen-multi-k}).
Indeed, if $k_i=i-1$, then $p^{i-1-k_i}=1$, so the congruence holds trivially, whereas  
if $k_i < i-1$, then $m_{i-1}=m_i$, and the congruence holds again. 

It remains to verify that $o_{p^i}(m_i)=(d_i)_{p'}$ 
for every integer $i$, $1 \le i \le e$. 
This is clear by the construction of $\vec{m}$ when 
$i=e$ or when $i < e$ and $k_{i+1}=i$. 
Assume that $i < e$ and $k_{i+1} < i$. 
Then $m_i=m_{i+1}$. Furthermore, we may assume w.l.o.g.~that 
$o_{p^{i+1}}(m_{i+1})=(d_{i+1})_{p'}$, and so we have that 
$$
o_{p^i}(m_i)=o_{p^i}(m_{i+1})=o_{p^{i+1}}(m_{i+1})=
(d_{i+1})_{p'},
$$
where the second equality follows by \eqref{eq:mod}. 
Finally, as $(d_1,\dots,d_t;\Pi)$ is the 
S-system corresponding to $\A_{G/E}$, 
$(d_{i+1})_{p'}=(d_i)_{p'}$, and so we conclude that $o_{p^i}(m_i)=(d_i)_{p'}$. 
\end{proof}

\begin{claim4} 
Let $X \in \S(\A)$ such that $X \not\subseteq H$ and $X \ne \{a,ab\}$. 
Then $X^{(m_i)}=Xb$, where $|\sg{X}|=2p^i$.
\end{claim4}
\begin{proof}

For an easier notation, we write $\bar{x}$ for the image of $x \in G$ under 
the canonical homomorphism $G \to G/E$, and let $\bar{Y}=\{\bar{x} : x \in Y \}$ 
for a subset $Y \subseteq G$. The set $\bar{X}$ is a basic set of the S-ring 
$\A_{G/E}=\S_{G/E}(d_1,\ldots,d_e;\Pi)$. According to \eqref{eq:Delta}, 
$$
\bar{X}=\bar{P}_i \setminus \bar{P}_j~\text{or}~\bar{X}=\orb_K(\bar{x}),
$$
where $0 \le j < i-1$, 
$x \in X$, and $K \le \aut(\bar{P_i})$ such that $|K|=d_i$. 
Observe that $\bar{X}^{(m_i)}=\bar{X}$. 
This is clear when $\bar{X}=\bar{P}_i \setminus \bar{P}_j$, and 
it follows from the condition $o_{p^i}(m_i)=(d_i)_{p'}$
when $\bar{x}=\orb_K(\bar{x})$. Then 
$(EX)^{(m_i)}=EX$, implying that $X^{(m_i)} \subseteq (X \cup Xb)$. 
Now, if $Xb=X$, then Claim~4 holds, thus we may assume that $Xb \ne X$. 
Since $\underline{\{b\}} \in \A$, $Xb$ is also 
a basic set of $\A$, and so we have $Xb \cap X=\emptyset$. 
On the other hand, $X^{(m_i)}$ is a basic set of $\A$ too, see 
Theorem~\ref{S-multi}(1), and therefore, it is sufficient to 
show that $X^{(m_i)} \ne X$. 

As $X \cap Xb=\emptyset$, we can write 
$$
X=X_1a \cup X_2ab, 
$$
where $X_1$ and $X_2$ disjoint subsets contained in $P_i$. 
Then $\sg{X^{[2]}}=P_i$, hence 
$\underline{P}_i \in \A$. Observe that Lemma~\ref{calB} can be applied to the S-ring 
$\A_{EP_i}$. As a result, we obtain that $X_1$ and $X_2$ 
are the basic sets of an S-ring, say $\B$ over the group $EP_i$, and also that 
$|X_1|=|X_2|$. 

Assume for the moment that $\bar{X}=\bar{P}_i \setminus \bar{P}_j$. 
If $X_1$ contains no element of order $2p^i$, then 
$$
|X_1| \le p^{i-1}-p^j < p^i-p^{i-1} \le |X_2|, 
$$
contradicting the condition that $|X_1|=|X_2|$. It follows that both $X_1$ and 
$X_2$ contain an element of order $2p^i$. Clearly, one of them must also 
contain an element of order $2p^{i-1}$, say $X_1$. Then applying 
Theorem~\ref{P2} to $\B_{G/E}$ and its basic set $\bar{X}_1$, we find that 
$P_i \setminus P_{i-1} \subseteq \bar{X}_1$.  Consequently,   
$X_2$ contains no elements of order $2p^i$, which is impossible. 

Thus $\bar{X}=\orb_{K}(\bar{x})$. On the other hand 
$\bar{X}_1$ and $\bar{X}_2$ are basic sets of $\B_{G/E}$ having the same size $|K|/2$.
We obtain that $\bar{X_i}=\orb_{\hat{K}}(\bar{y}_i)$, where $\hat{K} < K$, 
$|\hat{K}|=|K|/2$ and $y_1 \in X_1a$ and $y_2\in X_2ab$. 
Put $z=y_1^{m_i}$. Since $o_{p^i}(m_i)=(d_i)_{p'}$, it follows that 
$\bar{z} \in \orb_{\hat{K}}(\bar{y_2})$. This means 
$y_1^{m_i} \in X_2a \cup X_2ab$. On the other hand, 
as $y_1=x_1a$, this shows that $y_1^{m_i} \in X_2a$, and so 
$y_1^{m_i} \notin X$. We proved that $X^{(m_i)} \ne X$.
\end{proof}

Let $\varphi_{\vec{m}}$ be the permutation of $H$ defined in \eqref{eq:f}. 
Choose an arbitrary element $x \in S$ and let $X$ be the basic set of $\A$ 
containing $xa$. Suppose that $o(x)=p^i$ or $2p^i$. 
If $o(x)=2$, then $X=\{a,ab\}$. As $\underline{Sa} \in \A$, 
$X \subseteq Sa$, contradicting that $1 \notin S$. Thus $1 < i \le e$. 
By Claim~2, $P_{k_i}xa \subseteq X$. The coset 
$P_{k_i}x$ is a class of $\Pi(\k)$, hence by 
Proposition~\ref{M1}, $x^{\varphi_{\vec{m}}} \in (P_{k_i}x)^{(m_i)}$. 
Finally, by Claim~4, $X^{(m_i)}=Xb$. Combining all these together yields
$$
x^{\varphi_{\vec{m}}}a \in (P_{k_i}x)^{(m_i)}a=(P_{k_i}xa)^{(m_i)} 
\subseteq X^{(m_i)}=Xb \subseteq Sab.
$$
This shows that  $x^{\varphi_{\vec{m}}} \in Sb$, so 
$S^{\varphi_{\vec{m}}}=Sb$, as required.
\end{proof}
%--------------------------------------------------------------------------------------------------------------%


\begin{thebibliography}{mmm}
\bibitem{EP02}
S. Evdokimov and I. Ponomarenko,
On a family of Schur rings over a finite cyclic group,
St. Petersburg Math. J. 13 (2002) 441--451.
%
\bibitem{EP12}
S. Evdokimov and I. Ponomarenko,
Schurity of S-rings over a cyclic groups and generalised wreath product of 
permutation groups, 
Algebra Anal. 24 (2012) 84--127.
%
\bibitem{FH}
B. Fernandes and A. Hujdurovi\'c, 
Canonical double covers of circulants,
J. Comb. Theory Ser. B 154 (2022) 49--59.
%
\bibitem{HiMu}
M. Hirasaka and M. Muzychuk,
An elementary abelian group of rank $4$ is a CI-group,
J. Combin. Theory Ser. A 94 (2001) 339--362.
%
\bibitem{HuMi}
A. Hujdurovi\'c, \DH. Mitrovi\'c, 
Some conditions implying stability of graphs,
preprint arXiv:2210.15249 (2022).
%
\bibitem{HMW21a}
A. Hujdurovi\'c, \DH. Mitrovi\'c, and D. Witte Morris, 
On the automorphisms of the double cover of a circulant graph, 
Electron. J. Combin. 28 (2021) \#P4.43.
%
\bibitem{HMW21b}
A. Hujdurovi\'c, \DH. Mitrovi\'c, and D. Witte Morris, 
Automorphisms of the double cover of a circulant graph of valency at most $7$, 
to appear in Algebr. Comb.,  preprint https://arxiv.org/abs/2018.05164 (2021).
%
\bibitem{LM}
K. H. Leung and S. H. Man,
On Schur rings over cyclic groups II,
J. Algebra 183 (1996) 273--285.
%
\bibitem{MSZ}
D. Maru\v{s}i\v{c}, R. Scapellato, and N. Zagaglia Salvi, 
A characterization of particular symmetric (0,1) matrices. 
Linear Algebra Appl. 119 (1989) 153--162.
%
\bibitem{M}
M.  Muzychuk,
A solution of the isomorphism problem for circulant graphs,
Proc. London Math. Soc. 88 (2004) 1--41.
%
\bibitem{MP}
M.  Muzychuk, I. Ponomarenko
Schur rings,
European J. Combin. 30 (2009) 1526--1539.
%
\bibitem{P}
R. P\"oschel,
Untersuchungen von S-Ringen, insbesondere im Gruppenring von p-Gruppen,
Math. Nachr. 60 (1974) 1--27.   
%
\bibitem{QXZ19}
Y. L. Qin, B. Xia, S. Zhou, 
Stability of circulant graphs,
J. Combin. Theory Ser. B 136 (2019) 154--169.
%
\bibitem{QXZ21}
Y. L. Qin, B. Xia, S. Zhou, 
Canonical double covers of generalized Petersen graphs, and double 
generalized Petersen graphs
J. Graph Theory 97 (2021) 70--81.
%
\bibitem{Sch}
I. Schur,
Zur Theorie der einfach transitiven Permutationgruppen,
S.-B. Preuss. Akad. Wiss. Phys.-Math. Kl. 18 (1933) 598--623.
%
\bibitem{SM}
G. Somlai and M. Muzychuk,
The Cayley isomorphism property for $\Z_p^3 \times \Z_q$,
Algebr. Comb. 4 (2021) 289--299.
%
\bibitem{S}
D. B. Surowski, 
Stability of arc-transitive graphs, 
J. Graph Theory 38 (2001) 95--110.
%
\bibitem{T}
O. Tamaschke,
A generalisation of conjugacy in groups, 
Rend. Sem. Mat. Univ. Padova 40 (1968) 408--427.   
%
\bibitem{Wbook}
H. Wielandt, 
Finite permutation groups, Academic Press, New York 1964.
% 
\bibitem{W}
S. Wilson,
Unexpected symmetries in unstable graphs,
J. Combin. Theory Ser. B 98 (2008) 359--383.
%
\bibitem{W-M21}
D. Witte Morris, 
On automorphisms of direct products of Cayley graphs on abelian groups, 
Electron. J. Combin. 28 (2021) \#P3.5
%
\bibitem{W-M23}
D. Witte Morris, Automorphisms of the canonical double cover of a toroidal grid, preprint: https://arxiv.org/abs/2301.05396 (2023)
\end{thebibliography}
\end{document}